\documentclass[a4paper,12pt]{amsart}
\usepackage[margin=2.7cm]{geometry}
\usepackage{amsmath,amsfonts,amssymb,amsthm}
\usepackage{xcolor} 
\usepackage{shortmathj}

\newcommand{\CP}{\mathbb{CP}}

\newcommand{\Z}{\mathbb Z}

\newcommand{\X}{\widehat X}

\newcommand{\al}{\alpha}
\newcommand{\be}{\beta}

\newcommand{\La}{\Lambda}
\newcommand{\la}{\lambda}

\newcommand{\x}{\times}
 
\newcommand{\hra}{\hookrightarrow}
\newcommand{\del}{\partial}
\newcommand{\Tor}{\rm Tor}
\newcommand{\Spin}{\rm Spin}
\newcommand{\PD}{{\rm PD}}
\newcommand{\Hom}{\rm Hom}
\newcommand{\co}{\thinspace\colon}
\newcommand{\lra}{\longrightarrow}
\renewcommand{\int}{\rm int}
\newcommand{\PL}[1]{\textcolor{black}{#1}}

\newtheorem{thm}{Theorem}[section]
\newtheorem{lemma}[thm]{Lemma}
\newtheorem{cor}[thm]{Corollary}
\newtheorem{prop}[thm]{Proposition}
\newtheorem{question}[thm]{Question}

\theoremstyle{definition}
\newtheorem{defn}[thm]{Definition}
\newtheorem{rem}[thm]{Remark}
\newtheorem{rems}[thm]{Remarks}

\newtheorem*{ack}{Acknowledgments}

\begin{document}

\title{On almost complex embeddings of rational homology balls}

\author{Paolo Lisca}
\address{Dipartimento di Matematica, Largo Bruno Pontecorvo 5, 56127 Pisa, Italy} 
\email{paolo.lisca@unipi.it}	

\author{Andrea Parma}
\address{Dipartimento di Matematica, Largo Bruno Pontecorvo 5, 56127 Pisa, Italy} 
\email{andrea.parma94@gmail.com}	


\date{} 

\begin{abstract}
We use elementary arguments to prove that none of the Stein rational homology 4-balls shown by the authors and Brendan Owens to embed smoothly but not symplectically in the complex projective plane admit such almost complex embeddings. In particular, we show that those rational balls admit no symplectic embeddings in the complex projective plane without appealing to the work of Evans and Smith. 
\end{abstract}

\maketitle	

\section{Introduction}\label{s:intro} 

Let $p>q\geq 1$ be coprime integers and $B_{p,q}$ the Stein rational homology ball smoothing of the quotient singularity $\frac{1}{p^2}(pq-1,1)$. In~\cite{LP1} the authors extended work of Brendan Owens~\cite{Ow20} by exhibiting a subfamily 
\[
\{B(k,m),\ k\geq 0, m\geq 1\}\subset \{B_{p,q}\} 
\]
such that each $B(k,m)$ smoothly embeds in the complex projective plane. We later realized that the smooth embeddings constructed in~\cite{LP1} were obtained using certain special handlebody decompositions of $\CP^2$ \PL{called {\em horizontal}} and that every smooth, closed, orientable $4$-manifold admits \PL{horizontal} decompositions~\cite{LP2}. In~\cite{LP3} we use horizontal decompositions to prove the existence of many more smooth embeddings of the rational balls
$B_{p,q}$ into $\CP^2$. 

Work of Evans and Smith~\cite{ES18} -- based on Weimin Chen's adjunction formula for pseudoholomorphic curves in almost complex orbifolds~\cite{Ch04} -- implies that $B(k,m)$ admits no symplectic embedding in $\CP^2$. 
In fact, Evans and Smith show that $B_{p,q}$ embeds symplectically in $\CP^2$ if and only 
if there are integers $s$ and $t$ such that  
\begin{equation}\label{e:ES-condition}\tag{ES}
p^2+s^2+t^2 = 3pst\quad\text{and}\quad \pm q \equiv 3s/t \bmod p.
\end{equation}
Note that the above sign ambiguity is irrelevant because~\eqref{e:ES-condition} holds for $(p,q)$ if and only if it holds for $(p,p-q)$.  
In fact, $B_{p,q}$ and $B_{p,p-q}$ are symplectomorphic~\cite[Remark~2.8]{ES18}.  

The main purpose of this note is to show by elementary arguments 
independent of~\cite{ES18} that the smooth embeddings constructed in~\cite{Ow20,LP1} are 
not homotopic to almost complex embeddings. In particular, we deduce that the rational 
balls $B(k,m)$ admit no symplectic embedding in $\CP^2$ without appealing to~\cite{ES18}. 
We achieve this by associating to each collared, orientation-preserving 
topological embedding $j\co B_{p,q}\hra\CP^2$ an integer $c(j)$ which identifies the 
homotopy class of $j$ \PL{and a sign $h(j)\in\{\pm 1\}$ which, together with $c(j)$, determines whether the pull-back by $j$ of the standard complex structure of $\CP^2$ is homotopic to the Stein structure on $B_{p,q}$}. 
Recall that, given a closed, topological 4-manifold $X$, a topological 
embedding $j\co B_{p,q}\hra X$ is {\em collared} if the restriction $j|_{\del B_{p,q}}$ 
extends to a topological embedding $[-1,1]\x\del B_{p,q}\hra X$. The following is our main result. 

\begin{thm}\label{t:main}   
There is a collared, orientation-preserving topological embedding of $B_{p,q}$ into a closed, oriented topological 4-manifold $X$ homotopy equivalent to $\CP^2$ if and only if $pq-1$ is a quadratic residue modulo $p^2$. 
Moreover, to each collared topological embedding $j\co B_{p,q}\hra\CP^2$ one can associate an integer 
$0<c(j)<p^2$ such that $c(j)^2+1\equiv pq\bmod p^2$ \PL{and a sign $h(j)\in\{\pm 1\}$} having the following properties:
\begin{itemize} 
\item 
if $j_1, j_2\co B_{p,q}\hra\CP^2$ are two collared topological embeddings,   
then $j_1$ is homotopic to $j_2$ if and only if $c(j_1)\equiv c(j_2)\mod p$;
\item 
let $J_0$ denote the standard complex structure on $\CP^2$. Then, 
the pulled-back almost complex structure $j^*(J_0)$ is homotopic to the Stein structure 
on $B_{p,q}$ if and only if $\textcolor{black}{h(j)}q\equiv 3c(j)\bmod p$ and 
\begin{equation}\label{e:congr}
\dfrac{c(j)^2 (p^2-pq-1) - 1}{p^2}\equiv \frac{\textcolor{black}{h(j)}c(j)q + 3}{p}\bmod 2.
\end{equation} 
\end{itemize} 
\end{thm}

\begin{rems}\label{r:two}
(a) In view of Theorem~\ref{t:main} it is natural to ask whether 
there is a rational ball $B_{p,q}$ and two non-homotopic orientation-preserving 
embeddings $j_1, j_2\co B_{p,q}\hra\CP^2$. 
(b) Theorem~\ref{t:main} is consistent with the results of~\cite{ES18} in 
the following sense. For any symplectic embedding $j\co B_{p,q}\hra\CP^2$ the almost 
complex structure $j^*(J_0)$ is homotopic to the Stein structure 
on $B_{p,q}$ and therefore if~\eqref{e:ES-condition} holds 
the conditions of Theorem~\ref{t:main} must be satisfied by some integer $0<c<p^2$ \PL{and some sign $h\in\{\pm 1\}$}. 
Lemma~\ref{l:check} implies that this is indeed the case. 
\end{rems}

\begin{cor}\label{c:acemb}
Let $j\co B_{p,q}\hra \CP^2$ be a collared, orientation-preserving topological embedding. 
If $j^*(J_0)$ is homotopic to the Stein structure on $B_{p,q}$ then $q^2+9\equiv 0\bmod p$. 
In particular, $B(k,m)$ does not smoothly embed in $\ \CP^2$ as an almost complex manifold. 
A fortiori, $B(k,m)$ does not admit symplectic embeddings in $\ \CP^2$. 
\end{cor} 

\begin{proof} 
Let $j\co B_{p,q}\hra\CP^2$ be an embedding as in the statement 
such that $j^*(J_0)$ is homotopic to the Stein structure on $B_{p,q}$ and 
let $c(j)$ \PL{and $h(j)$ be the integers} provided by Theorem~\ref{t:main}. 
Since $c(j)^2+1\equiv pq\bmod p^2$, we have $c(j)^2\equiv -1\bmod p$, 
hence $\textcolor{black}{h(j)}q\equiv 3c(j)\bmod p$ implies $q^2\equiv -9\bmod p$. 
On the other hand, the rational balls $B(k,m)$ were 
shown~\cite{Ow20, LP1} to be of the form $B_{p,q}$ with $q^2+9$ not divisible 
by $p$ for each $k\geq 0$ and $m\geq 1$. The statement follows. 
\end{proof} 

\begin{rem}\label{r:open-stein}
By recent work of Gompf~\cite[Corollary~1.2]{Go20} the existence of a topological embedding $B_{p,q}\subset\CP^2$ implies that the (image of the) interior 
of $B_{p,q}$ is topologically isotopic to a Stein open subset $U\subset\CP^2$. By Corollary~\ref{c:acemb}, in the case of the smooth embeddings 
$B(k,m)\subset\CP^2$ of~\cite{Ow20, LP1} the Stein structure which exists on $\int(B(k,m))$ as a smoothing of a quotient singularity is not homotopic to the Stein structure pulled-back from $U$ by the time-1 map of the isotopy. 
\end{rem}

Earlier work of Gompf~\cite{Go13} implies that $B_{p,q}$ embeds holomorphically in $\CP^2$ 
if there is a smooth embedding $j\co B_{p,q}\hra\CP^2$ such that 
$j^*(J_0)$ is homotopic to the Stein structure on $B_{p,q}$. In fact, 
it follows by~\cite[Theorem~2.1]{Go13} that after a smooth ambient isotopy of $\CP^2$ the 
induced complex structure on the image of $j$ makes it a holomorphically embedded Stein handlebody. 
One can combine Theorem~\ref{t:main} with~\cite{Go13} to obtain the following. 

\begin{cor}\label{c:holom}
A smooth, orientation-preserving embedding $j\co B_{p,q}\hra\CP^2$ is homotopic 
to a holomorphic embedding if and only if $\textcolor{black}{h(j)} q\equiv 3c(j)\bmod p$ and Equation~\eqref{e:congr}
holds. 
\end{cor}

\begin{proof} 
We have observed above that by~\cite{Go13} if $j\co B_{p,q}\hra\CP^2$ is a smooth, 
orientation-preserving embedding, then $j$ is homotopic to a holomorphic embedding 
if and only if $j^*(J_0)$ is homotopic to the Stein structure on $B_{p,q}$. On the other 
hand, by Theorem~\ref{t:main} the latter condition on $j$ is equivalent 
to the stated congruences. 
\end{proof} 

\begin{rem}
We do not know which $B_{p,q}$'s smoothly embed in $\CP^2$, although many pairs $(p,q)$ are obstructed by Donaldson's Theorem~A~\cite{Do83}. Indeed, assuming $B_{p,q}\subset\CP^2$ one can construct a positive definite $4$-manifold of the form $W = P\cup\left(\CP^2\setminus B_{p,q}\right)$ for a suitable $4$-dimensional plumbing $P$ with positive definite intersection lattice $\La_P$. 
By Donaldson's theorem the intersection lattice $\La_W$ is standard, and arguing as in~\cite[Section~3]{Ow20} one can find a contradiction. 
On the other hand, it is not difficult to find $B_{p,q}$'s such that $p$ is not a Markov number, $p$ divides $q^2+9$ but Donaldon's theorem does not obstruct the existence of a smooth embedding. 
For instance, $B_{10,1}$ and $B_{73,8}$ are such rational balls. 
In fact, it is easy to check that \textcolor{black}{for $(c,h)\in\{(3,-1),(97,1)\}$ in the first case and $(c,h)\in\{(1998,1),(3331,-1)\}$ in the second case the conditions of Theorem~\ref{t:main} are satisfied, so that the balls $B_{10,1}$ and $B_{73,8}$ could conceivably admit holomorphic embeddings into $\CP^2$, although by~\cite{ES18} they admit no such symplectic embeddings.}
These examples have led us to the following question. 
\end{rem}

\begin{question}
Is there a $B_{p,q}$ which admits a holomorphic embedding but no symplectic embedding into $\CP^2$ ? 
\end{question}

\begin{ack}
The present work is part of MIUR-PRIN project 2017JZ2SW5. The authors wish to thank the referee for their accurate and helpful report. 
\end{ack}

\section{Embeddings into homotopy complex projective planes}\label{s:embhomotopy}

If $pq-1$ is a quadratic residue $\bmod\: p^2$ then the linking form of $-\del B_{p,q}$ 
is realized by the matrix $(p^2)$~\cite[Theorem~3.1]{Ed05}, so it follows from work of 
Boyer and Stong~\cite{Bo93, St93} that $-\del B_{p,q}$ is the boundary of a an oriented, 
simply connected topological 4-manifold $V$ with intersection form $(p^2)$. 
Gluing $B_{p,q}$ and $V$ along their boundaries produces a closed, oriented topological 
4-manifold $X$ homotopy equivalent to $\CP^2$ and containing a topologically 
embedded collared copy of $B_{p,q}$. This establishes one direction of the first sentence  
of Theorem~\ref{t:main}. To prove the other direction we will use the following   

\begin{prop}\label{p:generator} 
Let $X$ be a closed, oriented topological 4-manifold homotopy equivalent to $\CP^2$. 
Let $j\co B_{p,q}\hra X$ be a collared, orientation-preserving topological embedding  
and set $V := X\setminus\overline{j(B_{p,q})}$. Then, $H_1(V;\Z)=0$
and $H_2(V;\Z)\cong\Z$. Moreover, if $i\co V\hra X$ is the inclusion map, 
the subgroup $i_*(H_2(V;\Z))\subset H_2(X;\Z)\cong\Z$ has index $p$. 
\end{prop} 

\begin{proof}
A Mayer-Vietoris argument~\cite[Lemma~3.1]{Ow20} applied to the decomposition 
\[
\CP^2 = j(B_{p,q})\cup V
\]
gives $H_1(V;\Z)=0$ and $H_2(V;\Z)\cong\Z$. Let $g\in H_2(V;\Z)$ be a generator and 
\[
\al\in H_2(V,\del V;\Z)\cong H^2(V;\Z)\cong \Hom(H_2(V;\Z),\Z)
\]
a relative homology class such that $\al\cdot g=1$. Recall~\cite[\S~2.3]{ES18} 
that $H_1(B_{p,q})\cong \Z/p\Z$. Then, $p\del_* \al\in H_1(\del V;\Z)$ has 
zero image in $H_1(B_{p,q};\Z)$. This implies that $p\PD(\al)\in H^2(V;\Z)$ 
is the restriction of a class in $H^2(X;\Z)$, i.e.~$k$ times a generator $\La\in H^2(X;\Z)$ 
for some $k\in\Z$. Let $\ell := \PD(\La)\in H_2(X;\Z)$. Then,
\[
p=p\langle\PD(\al), g\rangle = k\langle\La,g\rangle = k(\ell\cdot i_* g). 
\]
Thus, we have $i_* g = d\ell$, where $d$ divides $p$. Exactness of the sequence   
\begin{align*}
H_2(\del V;\Z)=0\to H_2(V;\Z)\cong\Z \stackrel{(g\cdot g)\cdot}{\lra} H_2(V,\del V;\Z)\cong\Z\\ 
\to H_1(\del V;\Z)\cong \Z/p^2\Z \to H_1(V;\Z) = 0
\end{align*}
shows that $d^2=p^2$, and the statement follows.
\end{proof}

\begin{cor}\label{c:quadratic}
Let $X$ be a closed, oriented topological 4-manifold homotopy equivalent to $\CP^2$ 
and $j\co B_{p,q}\hra \CP^2$ a collared, orientation-preserving topological embedding. 
Then, $pq-1$ is a quadratic residue $\bmod\: p^2$. 
\end{cor}

\begin{proof}
Proposition~\ref{p:generator} implies that the intersection form of $V$ is represented by the matrix $(p^2)$, which therefore presents the linking form on $\del V = -\del B_{p,q} = L(p^2,p^2-pq+1)$. It easily follows (cf.~\cite[Theorem~3.1]{Ed05}) that $pq-1$ is a quadratic residue modulo $p^2$.  
\end{proof}

\section{An auxiliary 4-manifold and its intersection lattice}\label{s:X-hat}

In this section we establish some facts which will be used in Section~\ref{s:proof} to prove 
the second part of Theorem~\ref{t:main}. 
Let $j\co B_{p,q}\hra \CP^2$ be a collared, orientation-preserving topological embedding  
and let $V\subset\CP^2$ be the topological 4-manifold of Proposition~\ref{p:generator}. 
Let $R_{p,q}$ be the minimal resolution of the cyclic quotient singularity 
of type $\frac{1}{p^2}(pq-1,1)$. Note that there is a canonical identification $
\del R_{p,q} = \del B_{p,q}$ because $B_{p,q}$ is a smoothing of the same 
singularity. We use this identification to define the oriented topological 4-manifold 
\[
\X := R_{p,q}\cup V
\]
by gluing $\del R_{p,q}$ and $\del V$ along their boundaries. 
Let $n = b_2(R_{p,q})$. Using the Mayer-Vietoris sequence it is easy to check that $H_1(\X;\Z)=0$ 
and $b_2(\X) = n + 1$. By Poincar\'e duality and the Universal Coefficients Theorem we have 
\[
\Tor(H_2(\X;\Z))\cong\Tor(H_1(\X;\Z))=0. 
\]
Denote by $\La_R$, $\La_V$ and $\La_{\X}$, respectively, the free Abelian groups $H_2(R_{p,q};\Z)$, $H_2(V;\Z)$ and $H_2(\X;\Z)$ viewed as intersection lattices.
Recall that, as a smooth manifold, $R_{p,q}$ is the 4-dimensional plumbing of 2-disk bundles over spheres associated to a string of integers $(-a_1,\ldots,-a_n)$, where $a_i\geq 2$ for $i=1,\ldots, n$ and 
\begin{equation}\label{e:cfraction}
\textcolor{black}{[a_1,\ldots,a_n]:=\ }
a_1 - \cfrac{1}{a_2-\cfrac{1}{\cdots - \cfrac{1}{a_n}}}
= \frac{p^2}{pq-1}  
\end{equation}
The core 2-spheres of the plumbing $S_1,\ldots, S_n\subset R_{p,q}$ can be chosen to be smooth complex curves, and with their complex orientation they define the {\em vertex basis} 
\[
\{v_1=[S_1],\ldots,v_n=[S_n]\}\subset\La_R.
\] 
Let $D_n\subset R_{p,q}$ be a smoothly and properly embedded 
2-disk normal to $S_n$, oriented so that $S_n\cdot D_n=+1$. 
\textcolor{black}{Note that ``the last sphere'' $S_n$ is well-defined unless $n>1$ and 
\[
(a_1,a_2,\ldots, a_n)=(a_n,a_{n-1},\ldots, a_1), 
\]
which by~\cite[Lemmas~A.1 and~A.2]{OW77} holds if and only if $pq-1$ equals its inverse $\bmod\ p^2$, which is $p^2-pq-1$. Thus, $(a_1,\ldots, a_n)$ is palindromic if and only if $p=2q$. But, since $p$ and $q$ are coprime, this can happen only if $(p,q) = (2,1)$, and  
then $\dfrac{p^2}{pq-1} = 4 = [4]$, so $n=1$ and $S_1$ is well-defined.}

It is well-known and easy to check that the homology class $[\del D_n]$ is a generator of  
$H_1(\del R_{p,q};\Z)\cong\Z/p^2\Z$. Recall that by Proposition~\ref{p:generator} we have 
$H_2(V;\Z)\cong\Z$.  
\PL{
\begin{defn}\label{d:gen}
Let $g\in H_2(V;\Z)$ be the generator such that $i_* g\in H_2(\CP^2;\Z)$ is $p$ times the class of a complex line, where $i_*$ is the inclusion-induced map. 
Moreover, let 
\[
\al\in H_2(V,\del V;\Z)\cong H^2(V;\Z)\cong \Hom(H_2(V;\Z),\Z)\cong\Z
\]
be the (unique) relative homology class such that $\al\cdot g=1$. 
\end{defn}
}
Note that, since $g\cdot g = i_* g\cdot i_* g = p^2$ and we are assuming $p>1$, we must have $\del\al\neq 0\in H_1(V;\Z)=H_1(\del R_{p,q};\Z)$, otherwise $\al$ would be a multiple of $g$ and therefore $\al\cdot g\neq 1$.
\begin{defn}
Let $0< c(j)<p^2$ be the unique integer such that 
\[
c(j)[\del D_n] = \del\al\in H_1(\del R_{p,q};\Z).
\]
\end{defn}

\begin{lemma}\label{l:ext-basis}
There is an an element $v_{n+1}\in\La_{\X}$ such that 
$v_{n+1}\cdot i_* g = 1$ and 
\[
\{v_1,\ldots,v_n,v_{n+1}\}\subset\La_{\X}
\]
is a basis with associated Gram matrix $G = (v_i\cdot v_j)$ given by 
\[
G = 
\begin{pmatrix}
-a_1 & 1 & 0 & 0 & \cdots & 0 \\
1 & -a_2 & 1 & 0 & \cdots & 0\\
\vdots & \vdots & \vdots & \vdots & \vdots & \vdots \\
0 & \cdots & 1 & -a_{n-1} & 1 & 0\\
0 & \cdots & \cdots & 1 & -a_n & c(j) \\
0 & \cdots & \cdots & 0 & c(j) & -a_{n+1}
\end{pmatrix}.
\]
Moreover, $c(j)^2 \equiv pq-1\bmod p^2$ and $a_{n+1}  = \dfrac{c(j)^2 (p^2-pq-1) - 1}{p^2}$.
\end{lemma} 

\begin{proof}  
By construction, the pair $(c(j)\PD([D_n]),\PD(\al))\in H^2(R_{p,q};\Z)\oplus H^2(V;\Z)$ is 
mapped to zero by the difference of the restriction maps 
\[
H^2(R_{p,q};\Z)\oplus H^2(V;\Z)\to H^2(\del R_{p,q};\Z)\cong H^2(\del V;\Z).
\]
Therefore, by the cohomology Mayer-Vietoris sequence there is a homology 
class $v_{n+1}\in H_2(\X;\Z)$ such 
that $\PD(v_{n+1})\in H^2(\X;\Z)$ is sent to $c(j)\PD([D_n])$ 
by the restriction map $H^2(\X;\Z)\to H^2(R_{p,q};\Z)$
and to $\PD(\al)$ by the restriction map $H^2(\X;\Z)\to H^2(V;\Z)$. 

We claim that the classes $v_1,\ldots,v_n, v_{n+1}$ are a basis $\La_{\X}$. 
Since by construction $v_{n+1}\cdot g= 1$, the lattice $\langle g\rangle\subset H_2(V;\Z)$ generated by $g$ 
is {\em primitive} in $\La_{\X}$, i.e.~the quotient $\La_{\X}/\langle g\rangle$ is torsion-free. Since 
$g\cdot g=p^2$, the lattice $\langle g\rangle$ is nondegenerate. Applying~\cite[Lemma~A34]{Di12} we obtain 
\[
\det(\langle g\rangle^\perp) = \det(\langle g\rangle) = p^2. 
\]
Since $\langle v_1,\ldots, v_n\rangle$ is a finite-index sublattice of the non-degenerate lattice 
$\langle g\rangle^\perp$ having the same determinant, by~\cite[Lemma~A5]{Di12} we have 
$\langle v_1,\ldots, v_n\rangle = \langle g\rangle^\perp$. Thus, for each $\la\in\La_{\X}$ the class 
\[
\la - (\la\cdot g) v_{n+1} \in \langle g\rangle^\perp
\]
is a integral linear combination of $v_1,\ldots, v_n$. This shows  
that $\langle v_1,\ldots,v_n, v_{n+1}\rangle =\La_{\X}$, and since the classes $v_i$ 
are clearly independent the claim holds.  

Since $[D_n]\cdot v_n = 1$, by the definition of $v_{n+1}$ we have $v_{n+1}\cdot v_n = c(j)$. 
Setting $a_{n+1} := -v_{n+1}\cdot v_{n+1}$, we see that the Gram matrix $G = (v_i\cdot v_j)$ has 
the stated form. 
\PL{Using e.g.~\cite[Lemmas~A1 and~A2]{OW77} one can check that the numerator of the fraction $[a_1,\ldots, a_{n-1}]$ is the integer between $0$ and $p^2$ which is the inverse $\bmod\ p^2$ of $pq-1$, i.e. $p^2-pq-1$. 
} 
Moreover, since $\La_{\X}$ is unimodular of signature $(1,n)$ \PL{and the submatrix of $G$ given by the first $n-1$ rows and columns is negative definite}, we have  
\[
\det(G) = (-1)^n = -(-1)^n p^2 a_{n+1} - c(j)^2 (-1)^{n-1} (p^2-pq-1). 
\]
Thus, $c(j)^2 (-pq-1) \equiv 1\bmod p^2$, i.e.~$c(j)^2 \equiv pq-1\bmod p^2$, 
and the formula for $a_{n+1}$ follows. 
\end{proof} 

As before, let $g\in H_2(V;\Z)$ be the \PL{generator of Definition~\ref{d:gen}}. 
Then, by Lemma~\ref{l:ext-basis} we can write 
\[
i_* g = \sum_{i=1}^{n+1} b_i v_i \in \La_{\X}  
\]
for some $b_1,\ldots, b_{n+1}\in\Z$. 

\begin{lemma}\label{l:coordinates} 
There is a unique integer $h(j)\in\{\pm 1\}$ such that $b_1,\ldots, b_n$ are given by 
the recursive rule:
\[
b_1 = h(j)c(j),\quad b_2 = a_1 b_1, \quad 
b_s = a_{s-1} b_{s-1} - b_{s-2},\ s=3,\ldots, n.
\]
Moreover, $b_n = b_1 (p^2-pq-1)$ and $h(j)b_{n+1} = (a_n b_n - b_{n-1})/b_1 = p^2$. 
\end{lemma} 

\begin{proof}
Since $g\cdot g=p^2$, the lattice $\langle g\rangle$ is a finite-index sublattice of the non-degenerate lattice 
$\La_R^\perp=\langle v_1,\ldots, v_n\rangle^\perp$ 
having the same determinant, therefore by~\cite[Lemma~A5]{Di12} we have 
$\langle g\rangle = \La_R^\perp$. By Lemma~\ref{l:ext-basis},  
$\mathbf{b} := (b_1,\ldots, b_{n+1})^t\in\Z^{n+1}$ generates the kernel of 
the $n\x (n+1)$ matrix 
\[
M := 
\begin{pmatrix}
-a_1 & 1 & 0 & 0 & \cdots & 0 \\
1 & -a_2 & 1 & 0 & \cdots & 0\\
\vdots & \vdots & \vdots & \vdots & \vdots & \vdots \\
0 & \cdots & 1 & -a_{n-1} & 1 & 0\\
0 & \cdots & \cdots & 1 & -a_n & c(j) 
\end{pmatrix}
\]
viewed as a homomorphism $\Z^{n+1}\to\Z^n$. The system $M\mathbf{b} = 0$ 
consists of the $n$ equations  
\begin{equation}\label{e:eq1}
\begin{cases} 
	- a_1 b_1 + b_2 = 0,\\ 
	b_{s-2} - a_{s-1} b_{s-1} + b_s =0,\ s=3,\ldots, n,\\ 
	b_{n-1} - a_n b_n + c(j) b_{n+1} = 0.
\end{cases} 
\end{equation}
Since $\mathbf{b}$ is primitive the $b_i$'s are coprime, so the first $n-1$ equations of~\eqref{e:eq1} imply that 
$b_1$ divides $b_2,\ldots, b_n$. On the 
other hand, the last equation implies that $b_1$ divides $c(j) b_{n+1}$ and  
therefore, since the $b_i$'s are coprime, that $c(j) = h(j) b_1$ for some $h(j)\in\Z$.
Since $b_1$ divides $b_i$ for $i=1,\ldots, n$ and $c(j)\neq 0$, we can define $d_1,\ldots, d_n$ 
by the equations $b_i = b_1 d_i$, $i=1,\ldots, n$. 
Hence $d_1=1$ and the $d_i$'s satisfy the $n$ equations       
\begin{equation}\label{e:eq2}
\begin{cases} 
	- a_1 d_1 + d_2 = 0,\\
	d_{s-2} - a_{s-1} d_{s-1} + d_s =0,\ s=3,\ldots, n,\\ 
	d_{n-1} - a_n d_n + h(j) b_{n+1} = 0,
\end{cases} 
\end{equation}
which we may write in the form 
\begin{equation}\label{e:eq3}
	\dfrac{d_{i+1}}{d_i} = [a_i,\ldots,a_1],\ i=1,\ldots, n-1,\quad
	\dfrac{h(j) b_{n+1}}{d_n} = [a_n,\ldots, a_1].
\end{equation}
Note that~\eqref{e:eq2} and $d_1=1$ imply that $h(j)$ and $d_n$ are coprime. 
Hence, it follows from~\eqref{e:cfraction} and the last equation of~\eqref{e:eq3} 
that $h(j)$ divides $p^2$. But $h(j)$ also divides $c(j)$, which is coprime with 
$p$ because by Lemma~\ref{l:ext-basis} we have $c(j)^2\equiv -1 \bmod p$, 
therefore $h(j)=\pm 1$. It is easy to check by induction that, since $a_i>1$ for each $i$, we have  
\[ 
1 = d_1 < \cdots < d_n, 
\]
and the last equation of~\eqref{e:eq2} yields $h(j)b_{n+1}=(a_n b_n - b_{n-1})/b_1$. 
Finally, the last equality of~\eqref{e:eq3} implies that $h(j)b_{n+1} = p^2$ 
and by e.g.~\cite[Appendix]{OW77} $d_n (pq-1) \equiv 1\bmod p^2$, 
therefore $d_n = p^2 - pq - 1$. This concludes the proof. 
\end{proof}

\section{Embeddings into the complex projective plane}\label{s:proof} 

We continue using the notation of the previous sections. 

\begin{lemma}\label{l:homotopy}
Let $j_1, j_2\co B_{p,q}\hra\CP^2$ two collared topological embeddings.   
Then, $j_1$ is homotopic to $j_2$ if and only if $c(j_1) \equiv c(j_2)\bmod p$.
\end{lemma}

\begin{proof}
Recall that $\CP^\infty$ is a $K(\Z,2)$ and observe that, since $B_{p,q}$  
is homotopy equivalent to a $2$-dimensional 
CW-complex, the set $[B_{p,q},\CP^\infty]=H^2(B_{p,q};\Z)$ 
of homotopy classes of maps $B_{p,q}\to\CP^\infty$ 
is in $1-1$-correspondence with $[B_{p,q},\CP^2]$, the correspondence being given 
by composing a map $B_{p,q}\to\CP^2$ with the inclusion $\CP^2\subset\CP^\infty$. 
Hence the homotopy class of a map $j\co B_{p,q}\to\CP^2$
is determined by the pull-back $j^*(\PD(\ell))$, where $\ell\in H_2(\CP^2;\Z)$ is 
the class of a complex line. Therefore, $j_1$ is homotopic to $j_2$ if and only if 
$j_1^*(\PD(\ell)) = j_2^*(\PD(\ell))$. On the other hand, the cohomology exact sequence of the pair 
$(B_{p,q},\del B_{p,q})$ shows that the inclusion-induced map 
$H^2(B_{p,q};\Z)\to H^2(\del B_{p,q};\Z)$ is injective.  
Therefore $j_1^*(\PD(\ell)) = j_2^*(\PD(\ell))$ if and only if 
\[
j_1^* \PD(\ell)|_{\del B_{p,q}}
= j_2^*\PD(\ell)|_{\del B_{p,q}}. 
\]
Observe that, for each $k = 1,2$, we have 
$\PD(\ell)|_{\del V_k} = p\PD(\al)|_{\del V_k}$, 
where $V_k:=\overline{\CP^2\setminus j_k(B_{p,q})}$. 
Since by definition of $c(j_k)$ 
\[
j_k^*\PD(\ell)|_{\del V_k} =
p j_k^*\PD(\al))|_{\del V_k} = p c(j_k) \PD([\del D_n]),
\quad k=1,2,
\]
we conclude that $j_1$ is homotopic to $j_2$ if and only if $pc(j_1)\equiv pc(j_2)\bmod p^2$, 
i.e.~if and only if $c(j_1)\equiv c(j_2)\bmod p$. 
\end{proof}

\begin{lemma}\label{l:acstr}
Let $j\co B_{p,q}\hra\CP^2$ be a collared, orientation-preserving topological embedding. Then, 
the pulled-back almost complex structure $j^*(J_0)$ is homotopic to the Stein structure 
on $B_{p,q}$ if and only if $\textcolor{black}{h(j)}q\equiv 3c(j)\bmod p$ and Equation~\eqref{e:congr} holds. 
\end{lemma}

\begin{proof}
Let $s_0(\CP^2)\in\Spin^c(\CP^2)$ be the $\Spin^c$-structure associated to the standard complex structure on $\CP^2$, 
and $s_0(B_{p,q})\in\Spin^c(B_{p,q})$ the $\Spin^c$-structure associated to the Stein structure on $B_{p,q}$. 
Recall that, given any map $j\co B_{p,q}\to X$ there is an induced pull-back map 
$j^\sharp\co\Spin^c(X)\to\Spin^c(B_{p,q})$ between the sets of $\Spin^c$-structures on $X$ and $B_{p,q}$. 
Since $B_{p,q}$ is a $2$-complex, homotopy classes of almost complex structures on 
$B_{p,q}$ are in 1-1 correspondence with Spin$^c$ structures, with the correspondence given 
by sending an almost complex structure to the associated Spin$^c$ 
structure~\cite[Remark(a), p. 48]{Go97}. Therefore, $j^*(J_0)$ is homotopic to the Stein 
structure on $B_{p,q}$ if and only if $j^\sharp s_0(\CP^2) = s_0(B_{p,q})$. 
Therefore, to prove the theorem it suffices to show that $j^\sharp s_0(\CP^2) = s_0(B_{p,q})$ 
if and only if the stated congruences hold.

If $j^\sharp s_0(\CP^2) = s_0(B_{p,q})$ 
then, since $s_0(R_{p,q})|_{\del R_{p,q}} = s_0(B_{p,q})|_{\del B_{p,q}}$,  
$s_0(\CP^2)|_V$ extends to a $\Spin^c$-structure $\bar s_0\in\Spin^c(\widehat{X})$.
Let $\be = c_1(\bar s_0)\in H^2(\widehat{X};\Z)$. 
The class $\be$ restricts to $R_{p,q}$ as $c_1(R_{p,q})$, hence we have  
$\langle \be, i_* g\rangle = \langle c_1(\CP^2), i_* g\rangle = 3p$ 
because $c_1(\CP^2) = 3\PD(\ell)$, where $g\in H_2(V;\Z)$ is the generator of Proposition~\ref{p:generator}. Let $\{v_1^\#,\ldots, v_{n+1}^\#\}\subset H^2(\X;\Z)$ be the basis dual to $\{v_1,\ldots, v_{n+1}\}$. 
\PL{Since $v_1,\ldots,v_n\in H_2(R_{p,q};\Z)$ are homology classes of smooth complex curves of genus zero with their canonical orientation, by the classical adjunction formula we have}  
\[
\langle \be, v_i\rangle = \langle c_1(R_{p,q}), v_i\rangle = 2 - a_i\quad\text{for}\quad 
i=1,\ldots, n.  
\]
We can write 
\[
\be = \sum_{i=1}^n (2-a_i)v_i^\# + xv_{n+1}^\#\in H^2(\X;\Z)
\] 
for some $x\in\Z$ with $x = \langle \be, v_{n+1}\rangle\equiv a_{n+1}\bmod 2$ because 
$\be$ is characteristic. Since $i_* g=\sum_{i=1}^{n+1} b_i v_i$, we have 
\[
3p = \langle\be, i_* g\rangle = \sum_{i=1}^n (2-a_i)b_i + xb_{n+1}. 
\]
By Lemma~\ref{l:coordinates} 
\[
\sum_{i=2}^n b_i = \sum_{i=1}^{n-1} a_i b_i - \sum_{i=1}^{n-2} b_i 
\quad \Longrightarrow\quad  
\sum_{i=1}^n b_i = 
\sum_{i=1}^{n} a_i b_i - \sum_{i=1}^{n} b_i
 - b_1 pq,  
\]
where in the last equality we used that   
$b_1 p^2 = a_n b_n - b_{n-1}$ and $b_n = b_1(p^2-pq-1)$. 
Therefore 
\begin{equation}\label{e:betaval}
\sum_{i=1}^n(2-a_i) b_i = - b_1 p q,
\end{equation}
so by Lemma~\ref{l:coordinates} it follows that 
\textcolor{black}{
\[
3p = - b_1 p q + xb_{n+1} = h(j)p(-c(j)q+xp),
\]
}
therefore $\textcolor{black}{h(j)}q\equiv 3c(j)\bmod p$ and $x = \frac{c(j)q + 3\textcolor{black}{h(j)}}{p}$. 
Since $x = \langle \be, v_{n+1}\rangle$ 
must be congruent to $v_{n+1}\cdot v_{n+1} = -a_{n+1}$ $\bmod\: 2$ because $\be$ is characteristic, 
Equation~\eqref{e:congr} follows from Lemma~\ref{l:ext-basis}, and the first half of the proof is concluded. 

Conversely, let $j\co B_{p,q}\hra\CP^2$ be a collared, orientation-preserving topological embedding and 
suppose that $\textcolor{black}{h(j)} q\equiv 3c(j)\bmod p$ and~\eqref{e:congr} holds. Then, we can write 
\PL{$h(j)c(j)q + 3 = h(j)yp$}, where $y$ is an integer such that 
\[
y\equiv\dfrac{c(j)^2 (p^2-pq-1) - 1}{p^2}\bmod 2. 
\]
By Lemma~\ref{l:ext-basis} the class 
\[
\be := \sum_{i=1}^n (2-a_i)v_i^\# + yv_{n+1}^\#\in H^2(\X;\Z)
\]
satisfies $\langle\be,v_i\rangle\equiv v_i\cdot v_i\bmod\: 2$ for each $i=1,\ldots, n+1$, 
and therefore it is a characteristic class on $\widehat{X}$. Since $H_1(\widehat{X};\Z)=0$ we have  
$\be = c_1(s)$ for a unique $\Spin^c$ structure $s\in\Spin^c(\widehat{X})$,
and since $\be|_{R_{p,q}} = c_1(R_{p,q})$ the restriction  
$s|_{R_{p,q}}$ coincides with the $\Spin^c$-structure $s_0(R_{p,q})$ induced  
by the complex structure. Now observe that, by Lemma~\ref{l:coordinates} and our choice of $y$,  
\[
\langle\be, i_* g\rangle = \sum_{i=1}^n (2-a_i)b_i + yb_{n+1} 
= -b_1 pq + \textcolor{black}{yh(j)p^2 = p (-h(j)c(j)q + h(j)y p)} = 3p. 
\]

In view of Proposition~\ref{p:generator} this implies $\be|_V = c_1(\CP^2)|_V$, 
therefore, since $H_1(\CP^2;\Z)=0$, we have $s|V = s_0(\CP^2)|_V$. In 
particular, 
\[
s|_{\del V} = s_0(\CP^2)|_{\del V = j(\del B_{p,q})}.
\]
Hence, 
\[
j^\sharp s_0(\CP^2)|_{\del B_{p,q}} = s|_{\del R_{p,q}} = 
s_0(R_{p,q})|_{\del R_{p,q}} = s_0(B_{p,q})|_{\del B_{p,q}}.
\]
As observed in the proof of Lemma~\ref{l:homotopy}, 
since the map $H^2(B_{p,q};\Z)\to H^2(\del B_{p,q};\Z)$ induced by the 
inclusion $\del B_{p,q}\subset B_{p,q}$ is injective, 
so is the inclusion-induced map $\Spin^c(B_{p,q})\to\Spin^c(\del B_{p,q})$. 
Therefore 
\[
j^\sharp s_0(\CP^2)|_{\del B_{p,q}}=s_0(B_{p,q})|_{\del B_{p,q}} 
\quad\Longrightarrow\quad 
j^\sharp s_0(\CP^2)=s_0(B_{p,q}).
\]
This concludes the proof.
\end{proof} 

\begin{proof}[Proof of Theorem~\ref{t:main}]
The first sentence in the statement of Theorem~\ref{t:main} was established in Section~\ref{s:embhomotopy}. 
The second part of the statement follows from the combination of Lemmas~\ref{l:homotopy} and~\ref{l:acstr}.
\end{proof} 

We close the paper with the lemma referred to in Remark~\ref{r:two}(b).

\begin{lemma}\label{l:check}
If $p$ and $q$ satisfy Equation~\eqref{e:ES-condition} for some $s$ and $t$ then there 
exists an integer $c$ \PL{and a sign $h\in\{\pm 1\}$} satisfying the \PL{conditions} of Theorem~\ref{t:main}.
\end{lemma}

\begin{proof}
Given a triple $(p,s,t)$ satisfying~\eqref{e:ES-condition} we can choose an integer $c_0$ \PL{with $0<c_0<p^2$}, such that $c_0 \equiv\frac{s}{t}\bmod\: p$, therefore the condition $\textcolor{black}{h} q \equiv 3c \bmod p$ is satisfied \textcolor{black}{for some $h\in\{\pm 1\}$} by any $c$ of the form $c_0+kp$. 
It suffices to show that $c$ \PL{and $h$} satisfy the remaining congruences of Theorem~\ref{t:main} for some choice of $k$. Since $c_0^2 \equiv -1 \bmod p$ -- which implies that $c_0$ is coprime with $p$ -- we can write $c_0^2=ap-1$, so that $c^2+1 \equiv p(a+2c_0 k) \bmod p^2$. 
The right-hand side is congruent to $pq$ if and only if $2c_0k \equiv q-a \bmod p$, which becomes $2k \equiv (a-q)c_0 \bmod p$ after multiplying both sides by $c_0$. 
If $p$ is odd we can find $k$ by simply inverting $2$, while if $p$ is even there is a unique possibility for $k$ modulo $\frac{p}{2}$ and two possibilities modulo $p$ because $p\equiv 0\bmod 2$ implies that $(a-q)c_0$ is even. 
In fact, $q$ and $c_0$ are both odd because they are coprime with $p$. In particular, $c_0^2+1$ is congruent to $2\bmod 4$. Since $ap = c_0^2+1$, this shows 
that $a$ must be odd -- and $p\equiv 2\bmod 4$ -- so that $(a-q)c_0$ is even. 
We are left with verifying that $c$ satisfies \eqref{e:congr}. 
If $p$ is odd, replacing \PL{in~\eqref{e:congr}} all occurrencies of $p$ with $1$ we obtain the equivalent congruence \PL{$c^2 q + 1\equiv cq+1\bmod\: 2$}, which holds because $c^2 \equiv c \bmod 2$. 
If $p$ is even, by the argument given in the first part of the proof 
$\frac{p}{2}$ is odd. 
Moreover, replacing $c$ with $c+\frac{p^2}{2}$ changes the right-hand side $\bmod\: 2$ but not the left-hand side, so that exactly one between $c$ and $c+\frac{p^2}{2}$ satisfies the congruence. In fact, the right-hand side changes by \PL{$\dfrac{hpq}{2}$}, which is odd, while the left-hand side changes by 
\[
\dfrac{\left((c+\frac{p^2}{2})^2-c^2 \right)(p^2-pq-1)}{p^2}=
\left(c+\frac{p^2}{4}\right)(p^2-pq-1), 
\]
which is even because both $c$ and $\frac{p^2}{4}$ are odd.
\end{proof}

\bibliographystyle{amsplain} 
\bibliography{biblio}
\end{document}